\documentclass[centertags,leqno,11pt]{article}
\RequirePackage[mathscr]{eucal}
\RequirePackage{xspace,amssymb,mathrsfs,url,upref,verbatim}
\let\mathcal\mathscr

\usepackage{amsmath}
\usepackage{amssymb}
\usepackage{latexsym}
\usepackage{amsxtra}
\usepackage{amscd}
\usepackage{theorem}
\usepackage{fullpage}
\usepackage{titlesec}
\titleformat{\subsection}[runin]{\normalfont\bfseries}{\thesubsection.}{3pt}{}
\titleformat{\subsubsection}[runin]{\normalfont\bfseries}{\thesubsubsection.}{3pt}{}{}

\usepackage{graphics}
\usepackage{epic}

\theoremstyle{change}

{\theorembodyfont{\slshape}

\newtheorem{thm}{Theorem.}[section]
\newtheorem{cor}[thm]{Corollary.}

\newtheorem{prop}[thm]{Proposition.}

\newtheorem{defn}[thm]{Definition.}
}

{\theorembodyfont{\rmfamily}

\newtheorem{subsub}[thm]{}

\newtheorem{rem}[thm]{Remark.}

}

\newcommand{\proof}{\noindent {\bf Proof:\ }}
\newcommand{\Endproof}{\hspace*{\fill} $\Box$ \vspace{1ex} \noindent }

\makeatletter
\renewcommand{\subsection}{\@startsection{subsection}{2}%
{\z@}{-3.25ex plus -1ex minus-.2ex}{-1em}{\bf}} \makeatother

\newcommand{\PP}{\mathbb{P}}
\newcommand{\ZZ}{\mathbb{Z}}
\newcommand{\CC}{\mathbb{C}}

\newcommand{\QQ}{\mathbb{Q}}
\newcommand{\NN}{\mathbb{N}}
\newcommand{\FF}{\mathbb{F}}
\renewcommand{\AA}{\mathbb{A}}
\newcommand{\GG}{\mathbb{G}}

\newcommand{\OO}{\mathcal{O}}

\newcommand{\geom}{{\rm geo}}

\newcommand{\bQl}{{\overline{\QQ}_\ell}}

\newcommand{\x}{{\bf x}}

\newcommand{\GL}{{\rm GL}}
\newcommand{\SL}{{\rm SL}}

\newcommand{\Gal}{{\rm Gal}}
\newcommand{\Aut}{{\rm Aut}}
\newcommand{\Spec}{{\rm Spec\,}}

\newcommand{\Hom}{{\rm Hom}}

\newcommand{\MC}{{\rm MC}}
\newcommand{\SU}{{\rm SU}}

\newcommand{\cN}{{\mathcal N}}

\newcommand{\To}{\;\longrightarrow\;}

\newcommand{\Mapsto}{\;\longmapsto\;}

\newcommand{\diag}{{\rm diag}}

\newcommand{\pr}{{{\rm pr}}}

\newcommand{\J}{{\rm J}}

\newcommand{\Constr}{{\rm Constr}}

\newcommand{\Frob}{{\rm Frob}}

\newcommand{\SO}{{\rm SO}}

\renewcommand{\L}{{\mathcal L}}

\newcommand{\1}{{\bf 1}}
\newcommand{\cE}{{\mathcal E}}
\newcommand{\cF}{{\mathcal F}}
\newcommand{\cG}{{\mathcal G}}
\newcommand{\cH}{{\mathcal H}}

\newcommand{\cL}{{\mathcal L}}

\newcommand{\cO}{{\mathcal O}}

\newcommand{\rD}{{\rm D}}

\newcommand{\Ind}{{\rm Ind}}
\renewcommand{\cO}{{\mathcal{O}}}

\newcommand{\Gr}{{\rm Gr}}
\newcommand{\bT}{{\bf{T}}}

\newcommand{\cR}{{\cal R}}
\renewcommand{\Im}{{\rm Im}}
\newcommand{\smooth}{{\rm LocSys}}
\renewcommand{\Constr}{{\rm Constr}}

\numberwithin{equation}{section}
\numberwithin{thm}{subsection}

\theoremstyle{plain}

\title{ On Galois realizations of special linear groups}
\author{
   Michael Dettweiler and Stefan Reiter
  }

\begin{document}
\maketitle 
\begin{abstract} We study the determinant of certain \'etale sheaves 
constructed via  middle convolution 
in order to realize  special linear groups
regularly as Galois groups over $\QQ(t).$ \end{abstract}

\tableofcontents 


\section*{Introduction}\label{Introduction}

Recall that the { regular inverse Galois problem} is the following  question:\\

{\it \noindent Given a finite group $G,$ does  there exist a
Galois extension $L/\QQ(t)$ with $G\simeq \Gal(L/\QQ(t))$ such that additionally  $G\simeq \Gal(L/\overline{\QQ}(t))$ holds?}\\

\noindent If this condition holds for $G,$ then one says that $G$ {\it occurs regularly as Galois group over $\QQ(t).$}\footnote{It would be more accurate 
to say that the group $G$ occurs $\QQ$-regularly as Galois group over $\QQ(t)$.} 
The second  isomorphism, the {regularity condition,}  ensures that the field extension $L/\QQ(t)$ is geometric in the sense 
that it arises from 
a ramified  cover $f:X\to \PP^1_\QQ$ with  $\Aut(f)\simeq G.$ 

It
 follows from Hilbert's irreducibility theorem that a positive answer to 
the regular inverse Galois problem implies  a positive answer to the { inverse Galois problem}: {\it can every finite group $G$  be 
realized as Galois group of a Galois extension $L/\QQ$?} Both problems, however,  are  far from being solved, cf.~\cite{MalleMatzat}, \cite{Voelklein}. 
A weaker question, first posed by John Thompson, is the following: \\

\noindent {\it Given a finite field $\FF_q,$ is it  true that almost all finite groups of Lie type $G(\FF_q)$  occur regularly as Galois group over $\QQ(t)$?}\\

\noindent It follows from the work of V\"olklein and Thompson, V\"olklein (\cite{Voeklein92}, \cite{VoelkleinThompson}) and from our previous work \cite{DR99},\cite{DR00} that Thompson's question holds true under 
the further restriction to specific families of Lie type (like $\GL_n(\FF_q)$)  if $q$ is odd.    
It is the aim of this work to prove a similar result for the family of special linear groups (cf.~Thm.~\ref{thmrealierung1} and its corollary):\\

\noindent{\bf Theorem:} {\it Let $\FF_q$ be a finite field 
of odd order $q>3.$
Then the special linear group $\SL_{n}(\FF_q)$  
occurs regularly as Galois
group over $\QQ(t)$ if $n>8\varphi(q-1)+11,$ where $\varphi$ denotes Euler's $\varphi$-function.}\\

The proof relies on the Galois representations associated to certain
non-rigid \'etale sheaves  of rank two with finite monodromy.  Using two middle convolution
steps with quadratic Kummer sheaves, combined with a tensor operation with rank-one sheaves,  and using 
the permanence of having at most quadratic determinant under $\MC_{-\1}$ (\cite{De13},~Thm.~4.3.5), applied to these sheaves,
we obtain \'etale sheaves 
whose monodromy is contained in $\SL_n(\bQl).$ 
 The residual representations associated to these sheaves 
then give rise to the above result. \\

\thanks{We thank N. Katz for helpful remarks on an earlier version of this 
article.}

%

\section{Basic results and notation} 
Let in the following $R$ be a field or a normal integral domain which is finite  over $\ZZ,$ 
  let $X$ be 
 a connected, regular and separated scheme of finite type over $R,$ and let $\overline{x}$ be a geometric point of X. 

\subsection{Galois covers}\label{secbasics1}  (\cite{FreitagKiehl}, \cite{MalleMatzat}) 
 Any
finite  \'etale Galois cover $f:Y\to X$  with $G=\Aut(f)$ corresponds up to isomorphism to  a surjective homomorphism  of the \'etale fundamental group 
of $X$ onto $G$: $$\Pi_f:
\pi_1(X,\overline{x})\to G\leq {\rm Sym}(Y(\bar{x}))\quad {\rm with}\quad Y(\bar{x})=\Hom_X(\bar{x},Y).$$ 

\begin{subsub}{}\label{subsecgamma} Assume  that $R$ is a subring of $\CC$ and that $$X=\AA^1_R\setminus \x=\AA^1_R\setminus \{x_1,\ldots,x_r\}=\Spec(R[x][\frac{1}{(x-x_1)\cdots (x-x_r)}])$$ with  
$(x-x_1)\cdots (x-x_r)\in R[x]$ separable and $x_i\in \overline{{\rm Quot}(R)}$ \'etale over $\Spec(R).$  
Let  $$  \pi_1^{\rm top}(\PP^1(\CC)\setminus \{x_1,\ldots,x_r,x_{r+1}=\infty \})=\langle \gamma_1,\ldots,\gamma_{r+1}\mid \gamma_1\cdots \gamma_{r+1}=1\rangle
 \stackrel{\iota}{\To} \pi_1( \AA^1_{R}\setminus {\bf x})$$ be the natural inclusion, where 
 $\gamma_i\,(i=1,\ldots ,r+1)$ is a counterclockwise simple loop around $x_i$ as usual (cf.~\cite{MalleMatzat},~Chap.~I.1, \cite{FreitagKiehl}, Appendix~A). 
 The {\it monodromy tuple of $f:Y\to X$}  is by definition the tuple of elements
 $$ (\sigma_1,\ldots,\sigma_r,\sigma_\infty=\sigma_{r+1})\in G^{r+1}\quad \textrm{where}\quad \sigma_i:=\Pi_f(\iota(\gamma_i)),\quad  i=1,\ldots,r+1.$$
Note that by construction, the {\it product relation} $\sigma_1\cdots \sigma_{r+1}=1$ holds. 
Moreover, the operation of $G$  on $Y(\bar{x})$
is isomorphic to the regular representation of $G$ on itself. In this sense,  the homomorphism $\Pi_f$ will be viewed as homomorphism 
of $\pi_1(X,\overline{x})$ onto $G=\Aut(f).$ 
\end{subsub}

\begin{subsub}{}\label{subsubsec12} For $R$  finite  over $\ZZ,$ let $x\in X(\FF_{q^k}).$ Then the functoriality of $\pi_1$ yields a homomorphism$$ \pi_1(x,\overline{x})\simeq \Gal(\overline{\FF}_q/\FF_{q^k})\to \pi_1(X,\overline{x}) .$$ 
This leads to the notion of a (geometric) Frobenius element $\Frob_x$ in $\pi_1(X,\overline{x})$ by sending the profinite generator $\Frob_{q^k}$ 
of $\Gal(\overline{\FF}_q/\FF_{q^k})$ (which is 
 inverse to the arithmetic Frobenius $\overline{\FF}_{q^k}\to \overline{\FF}_{q^k},\, a \mapsto a^{q^k}$) to $\pi_1(X,\overline{x}).$
For $\tilde{x}$ another geometric point of $X$ there is an isomorphism $\pi_1(X,\overline{x})\stackrel{\beta_{\bar{x},\tilde{x}}}{\to} \pi_1(X,{\tilde{x}}),$ well 
defined up to inner automorphisms of $\pi_1(X,\overline{{x}})$. 
By sending $\Frob_x\in \pi_1(X,\overline{x})$ via this isomorphism to $\pi_1(X,{\tilde{x}})$ one obtains a Frobenius element, also denoted 
$\Frob_x,$ in the  group $\pi_1(X,{\tilde{x}}),$ well defined up to inner automorphisms. 
\end{subsub}

\subsection{Monodromy of \'etale sheaves.}\label{secfreitagk}   (\cite{DeligneWeil2}, \cite{FreitagKiehl}) 

\begin{subsub}{} Let in the following  $\cR$ be a ring used as coefficient field in \'etale cohomology
(like $\bQl,$ a finite extension of $\QQ_\ell,$  the valuation ring of such a field, or
the residue field of these).  Let $X$ be as in the previous section and let $\smooth(X,\cR),$ resp. 
$\Constr(X,\cR),$ denote the category of 
smooth (=lisse), resp. constructible, $\cR$-sheaves on $X.$ \\

For each geometric point $\bar{x}$ of $X,$ the association 
$L\in \smooth(X,\cR)\Mapsto L_{\bar{x}}$ establishes an equivalence of categories between $\smooth(X,\cR)$ and the category of finite continuous 
$\pi_1(X,\bar{x})$-modules. The monodromy representation of $L$ is by definition this representation $\rho_L:\pi_1(X,\bar{x})\to \Aut(L_{\bar{x}}).$ 
\end{subsub}

\begin{subsub}{} Let $f:Y\to X$ be  a Galois cover with associated homomorphism $\pi_f:
\pi_1(X,\bar{x})\to G=\Aut(f)$ as above and let  $V=\cR^n\, (n\in \NN).$ If 
$\rho:G \rightarrow \GL(V)$  is a representation, then one has a
sheaf $L=\L_{(f,\rho)}\in \smooth(X,\cR)$ associated  to the composition 
$$\rho_L=\rho\circ \Pi_f:\pi_1(X,\bar{x})\to \GL(V).$$  
\end{subsub}

\begin{subsub}{}\label{reminneroperation}   If   $y:\Spec(\FF_q)\to X$ is a closed point of $X,$ and 
if $L\in \Constr(X,\cR),$  then 
the stalk $L_{\bar{y}}$ is a $\pi_1(y,\bar{y})\simeq\Gal(\overline{\FF}_q/\FF_q)$-module in a natural way.
Hence
one has associated the characteristic polynomial 
$\det(1-\Frob_yt,L_{\bar{y}})$ to the Frobenius element $\Frob_y\in \pi_1(y,\bar{y}).$  Let $L\in \smooth(X,\cR)$ and let $\rho_L:\pi_1(X,\bar{x})\to \GL(V)$ be the monodromy 
representation of $L.$ By \cite{DeligneWeil2},~1.1.8, one has an equality of characteristic polynomials 
\begin{equation}\label{eqeqcharpol} \det(1-\Frob_yt,L_{\bar{y}})= 
\det(1-\rho_L(\Frob_y)t,V)\,,\end{equation} where on the right hand side, the Frobenius element $\Frob_y$ is viewed as an element (or rather a conjugacy class) 
in $\pi_1(X,\bar{x})$ via 
the isomorphism $\beta_{\bar{y},\bar{x}}$ from Section~\ref{subsubsec12}.
\end{subsub}

\begin{subsub}{} Let
$L\in \smooth(\AA^1_R\setminus {\bf x},\cR).$  
 Then  the {\it monodromy tuple}  of $L$ is defined as 
 $$\bT=\bT_L=(T_1,\ldots,T_{r+1})\in \GL_n(\cR)^{r+1},\quad T_i=\rho_L(\iota(\gamma_i))\, (i=1,\ldots,r+1),$$
 with $\iota:\pi_1(\AA^1(\CC)\setminus \x)\to \pi_1(\AA^1_R\setminus {\bf x})$ as in Section~\ref{subsecgamma}.
 \end{subsub}

 \subsection{Local monodromy}\label{seclocalmono2} (\cite{DeligneWeil2}, \cite{MalleMatzat})  Recall the notion of local monodromy: if $L$ is a smooth sheaf on an open subscheme $U$ of
 $X$ ($X$ a smooth and geometrically connected variety over a field $\kappa$) and if $x$ is a point 
 of $S=X\setminus U$ then the stalk $L_{\overline{\eta}_x}$ (with $\overline{\eta}_x$ denoting an algebraic closure of 
 the completion of the function field $\eta_x$ of $X$ w.r. to $x$) is a $\Gal(\overline{\eta}_x/{\eta}_x)$-module in a natural way: the {\it local monodromy 
 of $L$ at $x$}. The associated local monodromy representation is denoted $$\rho_{(x)}: \Gal(\overline{\eta}_x/{\eta}_x)\To\Aut(L_{\overline{\eta}_x})\simeq \GL_n(\cR).$$
 If $x$ is a closed point of $U$ then the stalk $L_{\overline{x}}$ is a $\Gal(\overline{\kappa}/\kappa(x))$-module. The associated 
 representation of $\Gal(\overline{\kappa}/\kappa(x))$ is denoted $\rho_x.$ 
 Note that for $X=\PP^1_\kappa,$ for $x\in \PP^1(\kappa),$ and for $L$ tame at $x,$ one has an isomorphism
 \begin{equation}\label{eqisoloc}
  \Gal(\overline{\eta}_x/{\eta}_x)^{\rm tame}=I_x^{\rm tame} \rtimes \Gal(\overline{\kappa}/\kappa)=\widehat{\ZZ}(1)(\overline{\kappa})\rtimes \Gal(\overline{\kappa}/\kappa)\,,
 \end{equation} where $\Gal(\overline{\eta}_x/{\eta}_x)^{\rm tame}$ denotes the tame quotient of $\Gal(\overline{\eta}_x/{\eta}_x)$
 and where
  $I^{\rm tame}_x$ denotes the tame inertia group at $x.$ 
If $x_i\in \PP^1(\CC)$ is as above then an 
image of a profinite generator $\gamma_i$ of $I_{x_i}^{\rm tame}$ in $\Aut(L_{\overline{\eta}_{x_i}})\simeq \GL_n(\cR)$ is conjugate 
to the $i$-th entry of the monodromy tuple $T_i$ (cf.~\cite{MalleMatzat}). Similarly as in Section~\ref{reminneroperation} 
one 
obtains a conjugacy class of morphisms $\Gal(\overline{\eta}_x/{\eta}_x)\to \pi_1(X,\tilde{x})$ (for $\tilde{x}$ another base point of $X$), describing 
the operation of $\Gal(\overline{\eta}_x/{\eta}_x)$ on $L_{\overline{\eta}_x}.$

\section{Construction of some smooth sheaves of rank $2$ with finite monodromy}
\subsection{The monodromy tuples} \label{subsec21}

In the following we use the following notation:
for $n\in \NN,$  let $(\zeta_n)_{n\in \NN}\in \overline{\QQ}$ be a system of primitive $n$-th roots of unity such that for $d\mid n$ 
one has $\zeta_d=\zeta_n^{n/d}.$ Let also $m\in \NN$ be a fixed integer $>2$  
and 
fix an embedding of $\bar{\QQ}$ into $\bQl.$\\

Let 
  $(T_1,\ldots,T_{r+1})\in \GL_2(\bQl)^{r+1}, r\geq 4,$ with
\begin{eqnarray*}
 T_i\,\,\,\,\,\, &=&\diag(\lambda_i,\lambda_i^{-1}),\quad i=1,\ldots,r-3, \,\textrm{with } 1\neq \lambda_i \in \bQl,\\
 T_{r-2}&=&\diag(1,-1) ,\\
 T_{r-1}&=&\left(\begin{array}{cc}
            0&1\\
           1&0 \end{array}\right),\\
 T_{r} \,\,\,\,\,\,&=&-(T_1 \cdots T_{r-1})^{-1},\\
 T_{r+1}&=& -\1_2,
\end{eqnarray*}
with $\1_n$ denoting the $n\times n$-identity matrix and with $\diag(\mu_1,\ldots,\mu_n)\in \GL_n(\bQl)$ denoting 
the diagonal matrix with diagonal entries $\mu_1,\ldots, \mu_n$ (in this order). \\

We assume in the following that the following 
conditions hold:
{\it \begin{enumerate}
 \item[(a)] $2\varphi(m)< r-4,$ 
 \item[(b)] the first $2\varphi(m)$
elements $\lambda_i$ run twice  through the primitive powers 
of $\zeta_m,$ and the remaining $\lambda_i$ are all equal to~$-1.$
\end{enumerate}}
 \noindent Under these conditions, we 
define $$\bT_{m,r}:=(T_1,\ldots,T_r,T_{r+1})\in\GL_2(\bQl)^{r+1} \quad \textrm{and}\quad Q_m:=\langle T_1,\ldots,T_{r+1} \rangle\leq \GL_2(\bQl).$$

\begin{rem}\label{remeigenspace2} {\it Note that then  the $r$-th component in $\bT_{m,r}$
is an element of order~$4,$ having a trivial $1$-eigenspace:
$$ T_r=\pm \left(\begin{array}{cc}
            0&-1\\
           1&0 \end{array}\right).$$ Note further 
           that the only components of $\bT_{m,r}$ with nontrivial invariants are the matrices 
           $T_{r-2}$ and~$T_{r-1}.$ }
      \end{rem}     

\subsection{Construction of the underlying sheaves}

\begin{subsub}{} 
It follows from the strong rigidity theorem (\cite{MalleMatzat},~Thm.~I.4.11) that there exists
an \'etale Galois cover $f:X\to \AA^1_\QQ\setminus {\bf \zeta}$ with ${\bf \zeta}=\{\zeta_m^d\mid d\in (\ZZ/m\ZZ)^* \}$ 
such that the monodromy tuple of $f$ is $$(\zeta_m,\zeta_m^{d_2},\ldots,\zeta_m^{d_{\varphi(m)}},1)\,\in (\bQl^\times)^{\varphi(m)+1},$$ 
with $ d_1=1,d_2,\ldots,d_{\varphi(m)}$ running through the elements of $(\ZZ/m\ZZ)^*,$ cf.~\cite{MalleMatzat},~Thm.~5.1.
The Galois cover $f,$ together with the embedding of $\mu_m$ into $\bQl^\times=\GL_1(\bQl),$ 
defines a smooth 
\'etale $\bQl$-sheaf $\cL_1$ on  $\AA^1_\QQ\setminus {\bf \zeta}$ of rank one. Let 
$${\bf \zeta'}=\{\pm \zeta_{2m}, \pm \zeta_{2m}^{d_2},\ldots,\pm \zeta_{2m}^{d_{\varphi(m)}}\}\,.$$ 
By pulling back $\cL_1$ along the map 
$\AA^1_\QQ\setminus {\bf \zeta'}\to \AA^1_\QQ\setminus \zeta,$ $x\mapsto x^2,$ one obtains a  smooth sheaf $\cL_2$ 
on $\AA^1_\QQ\setminus {\bf \zeta'}$ with monodromy tuple
 $$(\zeta_m,\zeta_m^{d_2},\ldots,\zeta_m^{d_{\varphi(m)}},
 \zeta_m,\zeta_m^{d_2},\ldots,\zeta_m^{d_{\varphi(m)}},1)\, \in (\bQl^\times)^{2\varphi(m)+1},$$ 
 up to a suitable renumeration of the elements in ${\bf \zeta'}.$ 
 \end{subsub}
 \begin{rem}\label{remisomstalks} {\it By construction, the sheaf $\cL_2$ has the property that
 for any $x\in \AA^1(\QQ)\setminus {\bf \zeta'},$ 
 under the chain of isomorphisms $\pi(x,\overline{x})\simeq \Gal(\overline{\QQ}/\QQ)\simeq \pi(-x,-\overline{x}),$ the 
$\Gal(\overline{\QQ}/\QQ)$-modules $(\cL_{2})_{\overline{x}}$ and $(\L_{2})_{-\overline{x}}$ are isomorphic. }
\end{rem}

\begin{subsub}{} {\it Let in the following  $\x=\{x_1,\ldots,x_{r-1}\}\in \overline{\QQ}^\times$ with $r-1>2\varphi(m)+2$ 
 be pairwise distinct points such that 
$$x_i=-x_{\varphi(m)+i}=\zeta_{2m}^{d_i}\quad (i=1,\ldots,\varphi(m), \, d_i\in (\ZZ/m\ZZ)^*\textrm{ as above} )$$ 
and such that for $i>2\varphi(m)$ the element $x_i$ is $\QQ$-rational and such that \begin{equation}\label{xr-1} x_{r-2}\neq -x_{r-1}.\end{equation}}

Since ${\bf \zeta'}\subset \x,$   we can view $\cL_2$ as a smooth sheaf on $\AA^1_\QQ\setminus \x$ by restriction. 
Using suitable 
quadratic covers, by the construction in Section~\ref{secbasics1}, there exist smooth sheaves $\cL_3,\cL_4$ on 
$\AA^1\setminus \x$ whose  monodromy tuples are $r$-tuples (the last component belonging to the point at $\infty$) of the form
$$ \bT_{\cL_3}=(1,\ldots,1,-1,\ldots,-1,1,1,\pm1)\quad {\rm and}\quad \bT_{\cL_4}=(1,\ldots,1,-1,-1,1),$$ 
resp., where in $\bT_{\cL_3}$ the entries $-1$ are at the positions $2\varphi(m)+1,\ldots,r-3.$ 
Let $\cL'=\cL_2\otimes \cL_3\otimes \cL_4$ and let $\cL''=(\cL_2\otimes \cL_3)^{\vee}$ be the dual of
$\cL_2\otimes \cL_3.$ 

Form the external tensor product $\cN=\cL'\boxtimes \cL''$ on $V=\AA^1_x\setminus \{x_1,\ldots,x_{r-1}\}\times 
\AA^1_y\setminus \{x_1,\ldots,x_{r-1}\}$ with respect to the canonical projections. Note that 
\begin{equation}\label{eqfundam} 
\pi_1(V,(\bar{x}_0,\bar{y}_0))=\pi_1(\AA^1_x\setminus \{x_1,\ldots,x_{r-1}\},\bar{x}_0)\times \pi_1(\AA^1_y\setminus \{x_1,\ldots,x_{r-1}\},\bar{y}_0)\end{equation} 
(where we view $\bar{x}_0, \bar{y}_0$ as complex points) comes equipped with the projections $\pi_x,\pi_y$ onto $\pi_1(\AA^1_x\setminus \{x_1,\ldots,x_{r-1}\},\bar{x}_0),$ resp. 
$\pi_1(\AA^1_y\setminus \{x_1,\ldots,x_{r-1}\},\bar{y}_0).$ 
In the following, we choose a base point $(\bar{x}_0,\bar{y}_0)$ of $V$ with $\QQ$-rational points 
$x_0,y_0$ satisfying $x_0\neq y_0.$ \\

Let $\gamma_{1,x},\ldots,\gamma_{r-1,x},\gamma_{\infty,x},$ resp. 
$\gamma_{1,y},\ldots,\gamma_{r-1,y},\gamma_{\infty,y},$ be standard counterclockwise generators of 
$\pi_1(\AA^1_x\setminus \{x_1,\ldots,x_{r-1}\}(\CC),\bar{x}_0),$ resp. 
$\pi_1(\AA^1_y\setminus \{x_1,\ldots,x_{r-1}\}(\CC),\bar{y}_0),$ viewed as elements in 
$\pi_1(\AA^1_x\setminus \{x_1,\ldots,x_{r-1}\},\bar{x}_0),$ resp. 
$\pi_1(\AA^1_y\setminus \{x_1,\ldots,x_{r-1}\},\bar{y}_0),$
as in Section~\ref{secbasics1}. Hence the  monodromy of $\cN$ is given by 
\begin{equation}\label{eqmonoN} \rho_\cN:\pi_1(V)\to \bQl^\times ,\quad \alpha\mapsto 
( \rho_{\cL_2\otimes \cL_3\otimes \cL_4})(\pi_x(\alpha))\cdot \rho_{\cL_2\otimes \cL_3}^{-1}(\pi_y(\alpha)).\end{equation} 
By Eq~\ref{eqfundam} we can view the elements $\gamma_{i,x},\gamma_{j,y}$ also as 
elements in $\pi_1(V,{(\bar{x}_0,\bar{y}_0)}).$
Hence \begin{equation}\label{eqMono4}
 \rho_\cN(\gamma_{i,x})= \rho_\cN(\gamma_{i+\varphi(m),x})=\zeta_m^{d_i}\,\, (i=1,\ldots,\varphi(m)),\quad 
 \rho_\cN(\gamma_{\varphi(m)+1,x})= \cdots=\rho_\cN(\gamma_{r-1,x})=-1,
\end{equation}
\begin{equation} \rho_\cN(\gamma_{i,y})= \rho_\cN(\gamma_{i+\varphi(m),y})=\zeta_m^{-d_i}\,\, (i=1,\ldots,\varphi(m)),\quad 
 \rho_\cN(\gamma_{2\varphi(m)+1,y})= \cdots=\rho_\cN(\gamma_{r-3,y})=-1\end{equation}
 and
 $$ \rho_\cN(\gamma_{r-2,y})=\rho_\cN(\gamma_{r-1,y})=1.$$ 
 \end{subsub}

\begin{subsub}{}  
Consider the canonical quotient map $$ h:\AA^2_{x,y}\to \AA^2_{s,t},\, (x,y)\mapsto (x+y,x\cdot y)$$ for the 
automorphism which switches the coordinates.
Under the map $h,$ the diagonal $\Delta\subset \AA^2_{x,y}$ is mapped to the 
conic $C:t-s^2/4=0,$ and a line $x-x_i=0$ (resp $y-x_i=0$) is mapped under $h$ to a tangent
$L_{x_i}:t-x_is+x_i^2=0$ to $C.$\\

Let $V'=V\setminus \Delta(V),$ with $\Delta(V)$ denoting the diagonal, let $\cN'=\cN|_{V'},$ and let 
$$W:=\AA^2_{s,t}\setminus( C\cup (\bigcup_{i=1,\ldots,r-1}L_{x_i})).$$ 
The map $h$ restricts to a quadratic 
\'etale cover 
$ \tau: V'\to W.$\\

The direct image $\cE=\tau_*(\cN')$ is a smooth rank-$2$ sheaf on $W$ whose monodromy representation is by construction the 
induced rank-$2$ representation $\rho_\cE=\Ind_{\pi_1(V')}^{\pi_1(W)}(\rho_\cN),$ where we view $\pi_1(V',(\bar{x}_0,\bar{y}_0))$ 
as a subgroup of $\pi_1(W,(\bar{s}_0,\bar{t}_0))$ (with 
$(\bar{s}_0,\bar{t}_0)=h(\bar{x}_0,\bar{y}_0)$).\\

Using a base point $(\bar{x}_0,\bar{y}_0)$ which is sufficiently close to the diagonal 
and by considering the punctured line $L$ through $(x_0,y_0)$ and $(y_0,x_0)$ one verifies the following: 
$$ \pi_1(V'(\CC),(\bar{x}_0,\bar{y}_0))=\langle \gamma_{i,x},\gamma_{i,y}, \gamma\mid i=1,\ldots, r-3\rangle, $$
where $\gamma$ is  a
path on $L(\CC)$ moving counterclockwise around $L\cap \Delta(\AA^2)(\CC)=(\frac{x_0+y_0}{2},\frac{x_0+y_0}{2})$ 
from $(\bar{x}_0,\bar{y}_0)$ to $(\bar{y}_0,\bar{x}_0)$ and back. \\

The image $\tau(L)$ is the parallel to the $t$-axis going through 
$(x_0+y_0,0).$ Let $\widetilde{\gamma}$ denote a simple loop in $\tau(L)$ around $\overline{\tau(L)}\cap C,$ 
represented by the non-closed half-twist in $L$  moving counterclockwise from $(x_0,y_0)$ to $(y_0,x_0)$ 
around $(\frac{x_0+y_0}{2},\frac{x_0+y_0}{2}),$ 
so that 
$\widetilde{\gamma}^2=\gamma.$ We can view $\pi_1(V')$ as a subgroup of $\pi_1(W)$ with $\gamma_{i,x}$ identified 
with a simple loop around the lines $L_{x_i},\, (i=1,\ldots,x_{r-1}).$ 
Then the 
fundamental group of $W$ is generated by  $\pi_1(V'),$ viewed as subgroup of $\pi_1(W),$ together with
 $\widetilde{\gamma}.$ By construction, 
$$ \gamma_{i,x}^{\widetilde{\gamma}}=\gamma_{i,y}\quad i=1,\ldots, r-1,$$ and vice versa.  
\end{subsub}
\begin{subsub}{} 
By the last remark we have 
\begin{equation}\label{eqmono2}
 \rho_\cE(\gamma_{i,x})=\Ind_{\pi_1(W)}^{\pi_1(V')}(\rho_\cN)(\gamma_{i,x})=\rho_\cN(\gamma_{i,x})\oplus \rho_\cN(\gamma_{i,x}^{\widetilde{\gamma}})
 =\rho_\cN(\gamma_{i,x})\oplus \rho_\cN(\gamma_{i,y}).
\end{equation}
With Eq.~\eqref{eqMono4} we obtain explicitly 
\begin{equation}\label{eqmono5}
 \rho_\cE(\gamma_{i,x})=\diag(\lambda_i,\lambda_i^{-1})\quad i=1,\ldots,r-3,
\end{equation}
with $\lambda_i=\lambda_{i+\varphi(m)}=\zeta_m^{d_i}\,(i=1,\ldots,\varphi(m))$ and $\lambda_i=-1$ for 
$i=2\varphi(m)+1,\ldots,r-3,$ 
and also we obtain 
\begin{equation}\label{eqmono6}
 \rho_\cE(\gamma_{i,x})=\diag(-1,1)\quad i=r-2,r-1.
\end{equation}
Since $\gamma_{i,x}^{\widetilde{\gamma}}=\gamma_{i,y}$ we conclude from Eq.~\eqref{eqmono2} that 
\begin{equation}\label{eqmono3}
 \rho_\cE(\widetilde{\gamma})=\left(\begin{array}{cc} 0&1\\
                       1&0
                      \end{array}\right).
\end{equation}
Let $\overline{Z}$ be the  connecting line in $\AA^2_{s,t}$  through $$z_r:=C\cap L_{x_{r-1}}$$ and through
$$z_{r-2}:=L_{x_{r-2}}\cap (\textrm{$t$-axis})=(0,-x_{r-2}^2).$$ Let $z_{r-1}$ denote the second intersection point of $\overline{Z}$ with 
$C.$ (Here we have used the condition in~\eqref{xr-1} to ensure that $\overline{Z}$ is not tangent to
$C.$)
 Let further $$Z=\overline{Z}\cap W\simeq \AA^1\setminus \{z_1,\ldots,z_r\}$$ with 
 $$z_i=\overline{z}\cap L_{x_i},\, i=1,\ldots,r-3.$$ 
 
Let $$\cF_{m,r}:=\cE|_{Z},$$ viewed as an object in $\smooth(\AA^1_\QQ\setminus \{z_1,\ldots,z_r\},\bQl).$ In case that $\cF_{m,r}$ has trivial local monodromy at $\infty$ we 
replace $\cF_{m,r}$ by a tensor product $\cF\otimes \cL_5,$ where $\cL_5$ is a rank-one sheaf
in $\smooth(  \AA^1_\QQ\setminus \{z_1,\ldots,z_r\},\bQl)$ whose monodromy tuple
is the $r+1$-tuple
$$ (1,\ldots,1,-1,-1).$$

\begin{prop}\label{proplocalmon3} There exist generators $\gamma_1,\ldots,\gamma_r$ of $\pi_1(\AA^1(\CC)\setminus \{z_1,\ldots,z_r\})$ such that 
 the monodromy tuple of $\cF_{m,r}$ with respect to these generators is the tuple 
   $\bT_{m,r}=(T_1,\ldots,T_{r+1})\in \GL_2(\bQl)^{r+1},$ specified in Section~\ref{subsec21}, assuming conditions a) and b). 
\end{prop}
\proof It follows from  the above description of the local monodromy of $\cE$ 
that the inertial local monodromy 
at the point $z_i$ is represented by $T_i$ ($i=1,\ldots,r+1$).  By \cite{DR00},~Lem.~7.2, the pure braid group acts transitively on the corresponding monodromy tuples, modulo diagonal 
conjugation with inner automorphisms from $Q_m.$ Hence, by an appropriate braiding, we can assume that 
there exist generators $\gamma_1,\ldots,\gamma_r$ of $\pi_1(Z)$ such that 
the monodromy tuple with respect to these generators coincides with~$\bT_{m,r}.$ 
 \Endproof

\begin{rem}\label{remF}{\rm  For any \'etale Galois cover $f:X\to \AA^1_\QQ\setminus {\bf x},$ by generic smoothness 
there 
exists an natural number $N$ and 
an \'etale Galois cover $$f_R: X_R\to \AA^1_R\setminus {\bf x}_R=\Spec(R[x][\frac{1}{(x-x_1)\cdots (x-x_r)}]\quad  (R=\ZZ[1/N])$$
such that $f$ is the base change of $f_R$ induced by the inclusion
$R \subseteq \QQ$ and such that the divisor $D$ associated to $(x-x_1)\cdots (x-x_r)$ 
is \'etale over the spectrum of $R.$ 

  Hence, the above sheaves $\cF_{m,r}$ extend to smooth sheaves on 
$\AA^1_R \setminus \{z_1,\ldots,z_r\},$  for $R=\ZZ[1/(N\cdot \ell)]$ with $N$ large enough, denoted  by the same symbols. 

With $D=\{z_1,\ldots,z_r\}\cup \infty$ 
and $j:\AA^1_R \setminus D \hookrightarrow \AA^1_R$ the inclusion
one sees that $j_*\cF_{m,r}[1]$ is an object in $\mathcal{T}(\AA^1_R,\bQl)_{R,D}$
in the sense of Def.~\ref{defT} of the Appendix to this article.}\end{rem}

Recall that the only components of the monodromy tuple $\bT_{m,r}$ with nontrivial invariants are the matrices 
           $T_{r-2}$ and~$T_{r-1},$ having a one-dimensional $1$-eigenspace.
            The operation of Frobenius elements on these 
            invariants is as follows:

\begin{prop}\label{propspecialize} Let  $\cF_{m,r}\in \smooth(\AA^1_R \setminus 
\{z_1,\ldots,z_r\},\bQl)$ where $R$ is as in Rem~\ref{remF}. 
Let $x,x'\in \AA^1_R(\FF_q)$ be $\FF_q$-points lying over $z_{r-1},$  $z_{r-2}$ (resp.),
where the characteristic of $\FF_q$ is also supposed to be $\neq \ell.$ 
  Then the elements
 ${\rm det}(\Frob_x,j_*\cF_{\overline{x}}) $ and  ${\rm det}(\Frob_{x'},j_*\cF_{\overline{x}'})$  are equal to  $\{\pm 1\}.$ \end{prop}

\proof To prove the claim for $\det(\Frob_x,j_*\cF_{\overline{x}}) $ it suffices to show that with  $h^{-1}(x)=(z,z),$ the stalk $j_*\cF_{\overline{x}}\simeq \cN_{(\overline{z},\overline{z})}$ 
is an at most  quadratic $\Frob_{(z,z)}$-module (hereby we can neglect the 
possible tensor product with $\cL_5$). 
 
By Eq.~\eqref{eqmonoN}, 
$$\rho_\cN(\Frob_{(z,z)})=\rho_{\cL_2}(\Frob_z)\rho_{\cL_3}(\Frob_z)\rho_{\cL_4}(\Frob_z)\cdot \rho_{\cL_2}^{-1}(\Frob_z)\rho_{\cL_3}^{-1}(\Frob_z)=\rho_{\cL_4}(\Frob_z)=\pm 1,$$
as claimed.

The stalk $j_*\cF_{\overline{x}'}$ is isomorphic to $\cN_{(-\overline{t}_0,\overline{t}_0)}$ for some $t_0\in \AA^1_R(\FF_q).$  Hence by  Eq.~\eqref{eqmonoN}, 
\begin{eqnarray}\nonumber 
\rho_\cN(\Frob_{(-t_0,t_0)})&=&
\rho_{\cL_2}(\Frob_{-t_0})\rho_{\cL_3}(\Frob_{-t_0})\rho_{\cL_4}(\Frob_{-t_0}) \rho_{\cL_2}^{-1}(\Frob_{t_0})
\rho_{\cL_3}^{-1}(\Frob_{t_0})\\
&=&
\rho_{\cL_2}(\Frob_{t_0})\rho_{\cL_3}(\Frob_{-t_0})\rho_{\cL_4}(\Frob_{-t_0})\rho_{\cL_2}^{-1}(\Frob_{t_0})
\rho_{\cL_3}^{-1}(\Frob_{t_0})\nonumber \\
&=&\rho_{\cL_3}(\Frob_{-t_0})\rho_{\cL_4}(\Frob_{-t_0}) 
\rho_{\cL_3}^{-1}(\Frob_{t_0})\nonumber\\
&=&\pm 1,\nonumber
\end{eqnarray}
where the second equality follows from the equality 
$$\rho_{\cL_2}(\Frob_{-t_0})=\rho_{\cL_2}(\Frob_{t_0})$$ holding by  the pullback construction of $\cL_2,$ cf.~Rem.~\ref{remisomstalks}. This proves the claim for 
$\Frob_{x'}.$
 \Endproof
\end{subsub}

\section{Galois realizations of special linear groups}

\subsection{Construction of the underlying sheaves via middle convolution}

Let $\cF=\cF_{m,r}$ be as in Rem.~\ref{remF}.
It follows from the existence of  suitable quadratic covers of $\AA^1_R\setminus \{z_1,\ldots,z_r\}$ (possibly by enlarging $R$)
that  
there exist smooth $\bQl$-sheaves  $\cN_1,\cdots,\cN_5$
on $\AA^1_R\setminus \{z_1,\ldots,z_r\}$
whose monodromy tuples are $r+1$-tuples of the form (resp.)
  $$
 \bT_{\cN_1}=(1,\ldots,1,1,-1,1,-1,1),,\quad \bT_{\cN_2}=\bT_{\cN_4}=(1,\ldots,1,1,-1,1,-1),$$
 $$ \bT_{\cN_3}=(1,\ldots,1,1,1,-1,-1,1,1,1),\quad  \bT_{\cN_5}=(1,\ldots,1,-1,1,1,-1,1).$$

Let $-\1:\pi_1(\GG_{m,R})\to \bQl^\times $ be the quadratic character
associated to the \'etale cover $\GG_{m,R}\to \GG_{m,R}, \,\,x\mapsto x^2,$ and to the
inclusion of $\Aut(f)\simeq\mu_2$ into $\bQl^\times.$ The latter data define a smooth sheaf $\cL_{-\1}$ on
$\GG_{m,R},$ which will be used in the following middle convolution steps
 (cf.~Rem.~\ref{remkummersh} of the Appendix).
In the following, we use the convolution $\MC_\chi$ of smooth sheaves as defined in Def.~\ref{defmcchiforsmooth} below.\\

We use the following notation: an expression like
$  (i,-i,J(2)^{2r-6},1),$
occurring  in Prop.~\ref{lema1} below,
denotes a matrix in Jordan canonical form having three Jordan blocks of length $1$ one for each eigenvalue $i,i^{-1},1$ (resp.)
and having $2r-6$ Jordan blocks of length~$2$ to the eigenvalue~$1,$ etc..\\

In the following, we use the middle convolution functor $\MC_\chi$ as defined
in Def.~\ref{defmcchiforsmooth} below.
The next result is an easy exercise using the numerology of
the middle convolution given in~\cite{Katz96},~Cor.~3.3.6 (cf.~\cite{DR07}, Prop.~1.2.1):

\begin{prop}\label{lema1} Let $m\in \NN_{>2}$ be even, let $r\in \NN$ with $2\varphi(m)\leq r-5,$ and let $\cF_{m,r}$ be the smooth sheaf on $\AA^1_R\setminus \{z_1,\ldots,z_r\}$
as in Rem.~\ref{remF}. Then the following holds:
\begin{enumerate}
\item The smooth $\bQl$-sheaf
\[ \cG_{1,m,r}:=   \cN_2\otimes \MC_{-\1}(\cN_1\otimes \MC_{-\1}(\cF_{m,r})) \]
has rank $n_1=4r-9.$ The Jordan form of the $i$-th entry $T_i$ of the monodromy tuple $\bT_{ \cG_{1,m,r}}$
is as follows (resp):
\begin{eqnarray*}
 {}(\lambda_i,\lambda_i^{-1},1^{4r-11})&& i=1,\ldots,r-3,\\
 (J(2)^{2r-6}, J(3)),&& i=r-2,\\
 {}(1,-1^{4r-10}),&& i=r-1,\\
 (i,-i,J(2)^{2r-6},1),&&i=r,\\
{}(1,\ldots,1),&& i=r+1.
\end{eqnarray*}
\item The smooth $\bQl$-sheaf
\[   \cG_{2,m,r}:= \cN_4\otimes \MC_{-1}(\cN_3\otimes \MC_{-1}(\cF_{m,r}\otimes \cN_5)) \]
has rank $n_2=4r-11.$ The Jordan form of the $i$-th entry $T_i$ of the monodromy tuple $\bT_{ \cG_{2,m,r}}$
is as follows (resp.):\begin{eqnarray*}
 {}(\lambda_i,\lambda_i^{-1},1^{4r-13}),&& i=1,\ldots,r-4,\\
  (\J(2)^{2r-6},1),&& i=r-3,\\
 (\J(3), \J(2)^{2r-8},1^2), &&i=r-2,\\
 {}(1,-1^{4r-12}), &&i=r-1,\\
 {}(i,-i,1^{4r-13}),&&i=r,\\
{}(1,\ldots,1)&&i=r+1\,.
\end{eqnarray*}
\item The smooth $\bQl$-sheaf
\[   \cG_{3,m,r}:= \MC_{-1}(\cN_5\otimes \MC_{-1}(\cF_{m,r})) \]
has rank $n_3=4r-10.$ The Jordan form of the $i$-th entry $T_i$ of the monodromy tuple $\bT_{ \cG_{3,m,r}}$
is as follows (resp.):\begin{eqnarray*}
 {}(\lambda_i,\lambda_i^{-1},1^{4r-12}),&& i=1,\ldots,r-4,\\
  (\J(3)^2,\J(2)^{2r-8}),&& i=r-3,\\
 {}(-1, 1^{4r-11}), &&i=r-2,r-1\\
  (i,-i,\J(2)^{2r-6}),&&i=r,\\
{}(-1,\ldots,-1)&&i=r+1\,.
\end{eqnarray*}
\item The smooth $\bQl$-sheaf
\[   \cG_{4,m,r}:=  \MC_{-1}(\cN_5\otimes \MC_{-1}(\cF_{m,r}\otimes \cN_5)) \]
has rank $n_4=4r-12.$ The Jordan form of the $i$-th entry $T_i$ of the monodromy tuple $\bT_{ \cG_{4,m,r}}$
is as follows (resp.):\begin{eqnarray*}
 {}(\lambda_i,\lambda_i^{-1},1^{4r-14}),&& i=1,\ldots,r-4,\\
  (\J(2)^{2r-6}),&& i=r-3,\\
{}(-1,1^{4r-13}), &&i=r-2,r-1,\\
 (i,-i,J(2)^{2r-8},1,1),&&i=r,\\
{}(-1,\ldots,-1)&&i=r+1\,.
\end{eqnarray*}
\end{enumerate}
\end{prop}

\subsection{Galois realizations of finite and profinite special linear groups}

  \begin{defn}{\rm 
  Let $H$ be a profinite group.
Then $H$ {\it occurs regularly
 as Galois group over $\QQ(t)$} if 
there exists a continuous 
  surjection $\kappa:\Gal(\overline{\QQ(t)}/\QQ(t))\to H $ 
such that  the restriction of $\kappa$ to $\Gal(\overline{\QQ(t)}/\overline{\QQ}(t))$ is surjective. }
\end{defn}

For an odd prime $\ell,$ let $q=\ell^k\, (k\in \NN_{>0}).$ 
Write $\cO_q$ for the valuation ring of the completion of $\QQ(\zeta_{q-1}),$ w.r. to a valuation $\lambda$ 
lying over $\ell.$ 

\begin{thm}\label{thmrealierung1} Let $\ell$ be an odd prime number and let $q=\ell^k\, (k\in \NN_{>0}).$ 
If $q>3$ then the special linear group $\SL_{n}(\cO_q)$  
occurs regularly as Galois
group over $\QQ(t)$ if $n>8\varphi(q-1)+11.$ \end{thm}

\proof By construction of the middle convolution for smooth sheaves 
(Def.~\ref{defmcchiforsmooth}), each  $\bQl$-sheaf $\cG_{i,q-1,r}\,(i=1,\ldots,4)$ 
of Prop.~\ref{lema1} is of the form $\widetilde{\cG}_{i,q-1,r}(-1)\otimes_{\cO_q}\bQl,$ 
where $\widetilde{\cG}_{i,q-1,r}$ is a smooth $\cO_q$-sheaf  (note the Tate twist by $-1$). Since $\cG_{i,q-1,r}$ is pure of weight~$2$ (since in each middle convolution step,
and on each $\FF_p$-fibre, the middle convolution 
is a higher direct image of an intermediate extension  which is pure of weight $0,$ resp. $1,$ cf.~\cite{DeligneWeil2},~\cite{BBD}),  the sheaf 
$\widetilde{\cG}_{i,q-1,r}$ 
is pure of  weight~$0.$ 

By Prop.~\ref{lema1}, the rank of $\cG_{i,q-1,r}\,(i=1,\ldots,4)$ is $n_1=4r-9,n_2=4r-11,n_3=4r-10,n_4=4r-12$ (resp.). 
We now divide the proof into the dimensions $n_1,n_2,n_3,n_4,$ beginning with $n_1:$

 Since $\cG_{1,q-1,r}$ is geometrically irreducible, the sheaf $\widetilde{\cG}_{1,q-1,r}$
is also irreducible and the monodromy tuple of $\widetilde{\cG}_{1,q-1,r}$ generates an (absolutely) irreducible subgroup of $\GL_{n_1}(\cO_q).$ 
It follows  from Prop.~\ref{lema1}~(i)
that the components of the monodromy tuple of $\widetilde{\cG}_{1,q-1,r}$ are contained in the special linear group. Hence  
the determinant $\det(\widetilde{\cG}_{1,q-1,r})=\Lambda^{n_1}(\widetilde{\cG}_{1,q-1,r})$  is a constant sheaf  of rank $1$ on $\AA^1_{R}\setminus \{z_1,\ldots z_r\}$
(with $R=\ZZ[\frac{1}{N}]$ for a large enough $N\in \NN$ as in Rem.~\ref{remF}).

It follows from Prop.~\ref{propspecialize} and  from absence of nontrivial invariants 
of the other local monodromies of $\cF_{q-1,r}$ 
 that  on each $\FF_p$-fibre $\AA^1_{\FF_p}\setminus \bar{\mathbf{z}} $ (where $p>N,p\neq \ell,$ and where $\bar{\mathbf{z}}$ denotes the reduction of the omitted divisor $\mathbf{z}:=\{{z}_1,\ldots,{z}_r\}$),
the  conditions of \cite{De13},~Thm.~{4.2.4}, are fulfilled for $F:=\cF_{q-1,r}|_{\AA^1_{\FF_p}\setminus \bar{\mathbf{z}}} :$
\begin{enumerate}
 \item The local geometric monodromy of $F$ at $\infty$ is scalar, given by the quadratic character $-\1:k^\times \to \overline{\mathbb{\QQ}}_\ell^\times,$
 but $F$ is not geometrically isomorphic to $\cL_{-\1}.$  
  \item  
   The $I_s^t$-module $\Gr^M(F_{\bar{\eta}_s})$ (= the semisimplification of the tame geometric inertia at $s,$ see~\cite{De13}, Section~3.2) 
   is self-dual for all $s$ in $\bar{\mathbf{z}}.$ 
   \item 
 For any  $x\in |\AA^1_{\FF_p}|$ 
   there exists an integer $m$ such that 
    $\det(\Frob_x,(j_*F)_{\overline{x}})=\pm q^m$
    (here: $m=0$).
 \end{enumerate} 
 It follows from \cite{De13},~Thm.~4.2.4, that these conditions 
 again hold for $\MC_{-\1}(\cF_{m,r})|_{\AA^1_{\FF_p}\setminus \bar{\mathbf{z}}} .$ This implies that these conditions are also  valid  for  $ \cN_1\otimes \MC_{-\1}(\cF_{m,r})|_{\AA^1_{\FF_p}\setminus \bar{\mathbf{z}}} ,$
since after tensoring with $\cN_1$ there are no inertial invariants at $\bar{z}_{r-2}$ and 
$\bar{z}_r.$ 
 Following  the construction process of $\cG_{1,q-1,r}$ via middle convolution
 (as in Prop.~\ref{lema1} (i)), applying $\MC_{-\1}$ and 
 \cite{De13},~Thm.~4.2.4 again (also noting that the tensor product with 
 $\cN_2$ does not change the property of having at most quadratic determinant up to Tate twist and also noting that  the underlying sheaf $\widetilde{\cG}_{1,q-1,r}$ has weight~$0$ on each $\FF_p$-fibre), 
 for each closed point $x\in |\AA^1_{\FF_p}\setminus \bar{\mathbf{z}}|$ one has   $$\det(\Frob_x,\widetilde{\cG}_{1,q-1,r}|_{\AA^1_{\FF_p}\setminus \bar{\mathbf{z}}})=\pm 1.$$
  Cebotarev's density theorem therefore  implies that  the determinant sheaf $\det(\widetilde{\cG}_{1,q-1,r})$ is the geometrically constant sheaf rank-one sheaf 
associated to an at most quadratic character $$\pi_1(\Spec(R),\Spec(\overline{\QQ}))\to \cO_q^\times.$$ 
Since the dimension $n_1=4r-9$ is odd, the  full arithmetic monodromy group of the sheaf $$\cH_{1,q-1,r}:=\widetilde{\cG}_{1,q-1,r}\otimes \det(\widetilde{\cG}_{1,q-1,r})$$ 
is hence contained in the group $\SL_{n_1}(\cO_q).$

Let $H^\geom=\Im(\rho^\geom_{\cH_{1,q-1,r}})\leq \SL_{1}(\cO_q)$ be 
the geometric monodromy group of $\cH_{1,q-1,r}$ and let 
$\overline{H}^\geom\leq \SL_{n_1}(\FF_{q})$ denote its image under the residual map on the coefficients (well defined up to semisimplification). 
 The middle convolution, as defined in  \cite{De13},  Def.~4.3.5, 
makes sense also over the coefficient field $\FF_q=\OO_q/\lambda$ and the basic properties (like preservation of irreduciblity and 
the effect on local monodromy) hold also in this case (for the irreduciblity one uses the same arguments as 
in \cite{De13},~Rem.~2.1.4, using$\mod\,\lambda$-coefficients, the effect of $\MC_\chi$ on the semisimplification of 
the$\mod \, \lambda$-local monodromy used below follows from the compatibility of the cohomological construction of $\MC_\chi$ 
with reduction$\mod\,\lambda$). Hence the group  $\overline{H}^\geom$ is an absolutely irreducible 
subgroup of $\SL_{n_1}(\FF_q),$ containing the negative of a reflection. Moreover, by \cite{DR99}, Prop.~6.6, 
$\overline{H}^\geom$ is primitive.

Hence, by the results of Wagner, Serezkin and Zalesskii (as collected in \cite{MalleMatzat},~Thm.~2.4), 
$\overline{H}^\geom$ contains a subgroup of type $\SU_{1}(\FF_{q'}),\SL_{n_1}(\FF_{q'})$ or the derived group $\Omega_{n_1}(\FF_{q'})$ of $\SO_{n_1}(\FF_{q'})$ 
(with $\FF_{q'}$ a subfield of $\FF_q$) as a normal subgroup. Note that the underlying dimension $n_1$ is $>8$ since $q>3,$ hence the exceptional cases 
in the list of Wagner, Serezkin and Zalesskii do not occur in our situation.  Note also that since the middle convolution $\MC_{-\1}$ preserves 
autoduality up to a Tate twist by Verdier  duality, we can exclude the groups $\Omega_{n_1}(\FF_{q'})$ 
since  the 
group $Q_m,$ viewed as a subgroup of $\GL_2(\FF_q)$ does not respect an orthogonal or a symplectic form. 
We can exclude the unitary groups $\SU_{n_1}(\FF_{q'})$ because they do not contain a bireflection of type 
$$T_1\mod \lambda=\diag(\zeta_{q-1},\zeta_{q-1}^{-1},1,\ldots,1) \mod \lambda.$$ 
Moreover, the Frobenius map $\Frob_{q'},$ for $q'$ a proper divisor of $q,$ does not stabilize the conjugacy class of the bireflection $T_1\mod \lambda.$
Therefore we have $q'=q$ and consequently  $\overline{H}^\geom= \SL_{n_1}(\FF_{q}).$ 
Since the residual map 
$\SL_{n_1}(\cO_q)\to  \SL_{n_1}(\FF_q)$ has the Frattini property (see \cite{Wei}, Cor.~A), we have 
$$H^\geom=\Im(\rho^\geom_{\cH_{1,q-1,r}})=\SL_{n_1}(\cO_q)=\Im(\rho_{\cH_{1,q-1,r}}),$$ where 
the last equality follows trivially from the inclusion of $\Im(\rho_{\cH_{1,q-1,r}})$ into $\SL_{n_1}(\cO_q).$ This proves  the 
claim for $n_1=4r-9$ since the absolute Galois group  $\Gal(\overline{\QQ(t)}/\QQ(t))$
surjects onto the \'etale fundamental group $\pi_1(\AA^1_{R}\setminus \{z_1,\ldots z_r\})$
appearing in the above monodromy representations.

The claim for $n_2$ follows from exactly the same arguments using the sheaf $\cG_{2,q-1,r}.$ 

The claim for $n_3,n_4$ uses the sheaves $\cG_{3,q-1,r} $ and $\cG_{4,q-1,r}$ and the same arguments to reduce to the case where the geometric and arithmetic monodromy group 
of the analogs $\widetilde{\cG}_{i,q-1,r}\in \smooth(\AA^1_R\setminus \{z_1,\ldots,z_r\}, \cO_q)  \, (i=3,4)$ of $\widetilde{\cG}_{1,q_1}$ are equal to the 
group $${\SL}^\pm_{n_i}(\cO_q)=\{A\in \GL_{n_i}(\cO_q)\mid \det(A)=\pm 1\}\quad i=3,4.$$ Note that ${\SL}^\pm_{n_i}(\cO_q)$
contains the special linear group $\SL_{n_i}\cO_q)$ as a subgroup 
of index~$2$ and that the only local monodromy matrices $T_i$ which do not lie in $\SL_{n_i}(\cO_q)$ are 
the elements $T_{r-2}$ and $T_{r-1},$ cf.~Prop.~\ref{lema1}. It follows therefore from the proof of \cite{MalleMatzat},~Thm.~I.5.3, applied successively to the tower 
of coverings belonging to $\rho_{\widetilde{\cG}_{i,q-1,r}}\otimes_{\cO_q} (\cO_q/\lambda^k),$
that the pullback $\widehat{\cG}_{i,q-1,r}$ 
of $\widetilde{\cG}_{i,q-1,r}$ to the quadratic cover 
$$\AA^1\setminus \x\to\AA^1_R\setminus \{z_1,\ldots,z_r\},\quad x\mapsto (x-z_{r-1})(x-z_{r-2}),$$  has geometric and arithmetic monodromy group equal to 
$\SL_{n_i}(\cO_q),$ proving the claim for $n_3$ and~$n_4.$  
\Endproof

\begin{cor} \label{corsldreia} Let $\FF_q$ be a finite field 
of odd order $q>3.$
Then the special linear group $\SL_{n}(\FF_q)$  
occurs regularly as Galois
group over $\QQ(t)$ if $n>8\varphi(q-1)+11.$ \Endproof
\end{cor}


\section{Appendix: Arithmetic middle convolution}\label{Appendix}

It is the aim of this section, which is basically a reformulation  of \cite{Katz96}, Chap.~4, to define an arithmetic version of the middle convolution 
which allows an application of the results of \cite{De13}  to our situation, where the omitted singularities are not contained in the ground field $\QQ.$

\begin{prop}\label{propsmoothj} Let $S$ be  an irreducible noetherian scheme, $X/S$ smooth,
and $D$ in $X$ a smooth $S$-divisor. For $F$  smooth on $X \setminus  D$ and tame along
$D, $ and for $j: X\setminus D\to X$ and $i: D \to X$  denoting the inclusions, the following holds:
\begin{enumerate}
\item formation of $j_*F$  and of $Rj_*F$  on $X$ commutes with arbitrary
change of base on $S,$
\item  the sheaf $i^*j_*F $ on $D$ is smooth, and formation of $i^*j_*F $ on $D$
commutes with arbitrary change of base on $S.$
\end{enumerate}
\end{prop}
\proof \cite{Katz96},~Lem.~4.3.8.\Endproof

Recall from \cite{KatzLaumon} that a scheme is called {\it good} if it admits a map 
of finite type to a base scheme $S={\rm Spec}(R)$  which is regular of dimension at most
one. For  good schemes X and $\ell$  a fixed prime number, invertible in $X,$
one has the triangulated category 
$\rD^b(X,\bQl)$ , which
admits the full Grothendieck formalism of the six operations (\cite{DeligneWeil2}, \cite{KatzLaumon}).

Let  $R$  be a normal noetherian integral
domain in which our fixed prime $\ell$  is invertible so that $S={\rm Spec}(R).$  Let 
$\AA^1_R=  \Spec(R[x])$ and let $D$ denote a smooth $S$-divisor defined by the 
vanishing of a separable monic polynomial $D(x)\in R[x]$ plus the divisor at~$\infty.$  

One says that an object $K\in \rD^b_c(\AA^1_R,\bQl)$ is {\it adapted to the stratification 
$(\AA^1\setminus D,D)$} if each of its cohomology
sheaves is smooth when restricted either to $\AA^1_R\setminus D$ or to $D$
(\cite{Katz96},~(4.1.2), \cite{KatzLaumon},~(3.0)).  

\begin{defn}\label{defT} Let $\mathcal{T}(\AA^1,\bQl)_{R,D}$ denote the category formed by the
 objects $K$ in $\rD^b_c(\AA^1_R,\bQl)$  of the form $j_*F[1],$ where 
 $j:\AA^1_R\setminus D\hookrightarrow \AA^1_R$ denotes the inclusion and 
 $F$ is smooth on $\AA^1_R\setminus D,$  such that the following holds: 
 
\begin{enumerate}
\item For  $k$ an algebraically closed field and $R\to k$ a ring homomorphism the restriction
$F|_{\AA^1_k\setminus D_k}$ is smooth, irreducible and nontrivial.
\item  
The sheaf $F|_{\AA^1_k\setminus D_k}$ has  at least three non-smooth points in $D_k$ (including $\infty$).
\end{enumerate}

Let $\mathcal{T}(\AA^1,\bQl)_R$ denote the category of sheaves $F$ on $\AA^1_R$ for which there exists 
a $D$ such that $F\in \mathcal{T}(\AA^1,\bQl)_{R,D}.$
\end{defn}

By the previous result, each $K\in \mathcal{T}(\AA^1,\bQl)_{D,R}$ is adapted to the  stratification 
$(\AA^1\setminus D,D).$ Moreover, the restriction of  $K\in \mathcal{T}(\AA^1,\bQl)_R$ to each geometric fiber 
$\AA^1_k$ is an intermediate
 extension of an irreducible smooth sheaf and is hence perverse
(cf. \cite{Katz96}, Chap.~4, and \cite{De13}, Section~1.2).\\

\begin{rem}\label{remkummersh} Let $N$ be a natural number $>1$ and let $R$ be as above such that $R$ contains 
a primitive $N$-th root of unity and such that $N$ is invertible in $R.$ 
Consider the \'etale cover 
$f:\GG_{m,R}\to \GG_{m,R}, \,x\mapsto x^N,$ with automorphism group $\mu_N$  and let $\chi: \mu_N\to \overline{\mathbb{\QQ}}_\ell^\times $ be a character.  
The latter data define a smooth sheaf $\cL_{\chi}$ on 
$\GG_{m,R},$ by pushing out the so obtained $\mu_N$-torsor by $\chi^{-1}.$

 Note that for the natural embedding
 $$-\1:\mu_2=\{\pm 1\}\hookrightarrow \overline{\mathbb{\QQ}}_\ell^\times$$ 
one obtains in this way a smooth sheaf $\cL_{-\1}$ on $\GG_{m,\ZZ[1/(N\cdot \ell)]}$ for any even $N.$ 
Then on each $\FF_q$-fibre ($q$ prime to $N\cdot \ell$), the restriction $\cL_{\chi}|_{\GG_{m,\FF_q}}$ is obtained by the same procedure
by first considering 
$f_{\FF_q}:\GG_{m,\FF_q}\to \GG_{m,\FF_q}, \,x\mapsto x^2,$ with automorphism group $\mu_2$  and by taking the 
same character $-\1: \mu_2\to \overline{\mathbb{\QQ}}_\ell^\times.$ 
By looking at Frobenius traces,  the sheaf $\cL_{-\1}|_{\GG_{m,\FF_q}}$ coincides with the usual Kummer sheaf associated 
to the quadratic character of $\mathbb{G}_m(\FF_q),$ see \cite{De13}, Section~1.4, and \cite{Laumon}. 
\end{rem}

Let $j:\AA^1_{R}\times \AA^1_{R}\hookrightarrow \PP^1_{R}\times \AA^1_{R}$
denote the inclusion and let $\overline{\pr}_2:\PP^1_{R}\times \AA^1_{R}\to \AA^1_{R}$ be the second projection. 

Following \cite{Katz96},  
for a nontrivial character $\chi$ as above, 
define the {\it middle convolution} of $K\in \mathcal{T}(\AA^1,\bQl)_R$ with $j'_*\cL_\chi[1]\,$  as follows (where 
$j'$ denotes the inclusion of $\GG_m$ into $\AA^1$ and where  $\tau_k$ denotes the natural truncation functor), cf.~\cite{Katz96} ~(4.3.2):
\begin{equation}\label{eqdefmc} \MC_\chi(K)=R\overline{\pr}_{2*}(\tau_{\leq -2}Rj_{*}(\pr_1^*K\boxtimes j'_*\cL_\chi(t-x)[1]))=R\overline{\pr}_{2*}(j_{*}(\pr_1^*K\boxtimes j'_*\cL_\chi(t-x)[1])),\end{equation}
where $\cL_\chi(t-x)$ denotes the pullback of $\cL_\chi$ along the map $t-x$
(here the second  equality holds by construction since, locally at the divisor at 
$\infty,$ the perverse sheaf $\pr_1^*K\boxtimes j'_*\cL_\chi(t-x)[1]$ is a sheaf placed
in cohomological degree $-2$).

%

\begin{thm}\label{thmrelMCchi} \begin{enumerate}
\item For $K\in \mathcal{T}(\AA^1,\bQl)_{R,D},$ the middle convolution
$\MC_\chi(K)$ is again an object of $ \mathcal{T}(\AA^1,\bQl)_{R,D}.$ 

\item Formation of $\MC_\chi$ commutes with arbitrary change of base. 
Especially, on each geometric fiber 
$\AA^1_k,$ with $k$ either an algebraically closed field, one has 
$$ \MC_\chi(K)|_{\AA^1_k}=\MC_\chi(K|_{ \AA^1_k}),$$
cf.~\cite{Katz96}, Prop.~2.9.2. \end{enumerate}
\end{thm}

\proof The second claim follows from the same arguments as in 
 \cite{Katz96}, (4.3.2)--(4.3.6). 
 The first claim follows using the same arguments as in the proof of \cite{Katz96}, Thm. 4.3.11.
 \Endproof
 

In view of the previous result, one can define $\MC_\chi$ also for 
 constructible and smooth sheaves:

\begin{defn}\label{defmcchiforsmooth}  Let $R,D,$ and $\chi$ be 
as above. \begin{enumerate}
\item Let $G$ be a constructible  $\bQl$-sheaf on $\AA^1_R$
such that $G[1]\in \mathcal{T}(\AA^1,\bQl)_{R,D}.$  Then the {\it middle convolution} of $G$ with respect to $\chi$
is defined as the constructible sheaf
\begin{equation}\label{eqdefmc2} \MC_\chi(G)=\MC_\chi(G[1])[-1]=\cH^{-1}(\MC_\chi(G[1]))\,.\end{equation} \label{defmc33} For $R=k$ an algebraically closed field this is Katz'
middle convolution functor ${\rm MC}_\chi,$ see \cite{Katz96}, (5.1.5).
\item For $F$ a smooth sheaf on $\AA^1_R\setminus D$ such that $j_*F[1]
\in \mathcal{T}(\AA^1,\bQl)_{R,D}$ define then $\MC_\chi(F)$ to be the smooth sheaf 
\begin{equation}\label{eqdefmc4} \MC_\chi(F)=\MC_\chi(j_*F)|_{\AA^1_R\setminus D}\,.\end{equation}
\end{enumerate}
\end{defn}
%
%
%

    \bibliographystyle{plain}
                    \bibliography{p1}

\end{document}